\documentclass[a4paper,12pt]{article}
\include{epsf}
\usepackage{latexsym}
\usepackage{amssymb}
\usepackage{amsmath}
\usepackage[dvips]{graphicx}
\newlength{\cqfd}
\def\R{\mathbb{R}}

\def\N{\mathbb{N}}

\newtheorem{proposition}{Proposition}

\newcommand{\bnu}{\overline{\nu}}
\newcommand{\bd}{\overline{\delta}}
\newcommand{\bq}{\overline{q}}

\newcommand{\LL}{\mathcal{L}}
\newcommand{\LLkdv}{\tilde{\mathcal{L}}}
\newcommand{\LLs}{\mathcal{L}_\star}
\newcommand{\LLkdvs}{\tilde{\mathcal{L}}_\star}

\newcommand{\ks}{k_\star}
\newcommand{\bqs}{\bq_\star}

\newcommand{\Us}{U_\star}
\newcommand{\Ms}{M_\star}

\newcommand{\ps}{p_\star}

\newcommand{\cs}{c_\star}

\title{\bf Modulated wave trains in generalized Kuramoto-Sivashinksi equations}

\begin{document}
 \maketitle 
 \begin{center}
{\large Pascal Noble \footnote{ Universit\'e de Lyon, Universit\'e Lyon 1
Institut Camille Jordan, UMR CNRS 5208
 43, blvd du 11 novembre 1918,
F - 69622 Villeurbanne Cedex, France; noble@math.univ-lyon1.fr:
Research of P.N. was partially supported by French ANR project
no. ANR-09-JCJC-0103-01}  L.Miguel Rodrigues  \footnote{Universit\'e de Lyon, Universit\'e Lyon 1
Institut Camille Jordan, UMR CNRS 5208
 43, blvd du 11 novembre 1918,
F - 69622 Villeurbanne Cedex, France; rodrigues@math.univ-lyon1.fr}}
\end{center}
\vspace{0.2cm}
 \begin{center}
\end{center}

\begin{center}
{\bf 2000 MR Subject Classification}: 35B35.
\end{center}

\vspace {0.5cm}
\begin{center}
{\bf Abstract.}
\end{center}
This paper is concerned with the stability of periodic wave trains in a generalized Kuramoto-Sivashinski (gKS) equation. This equation is useful to describe the weak instability of low frequency perturbations for thin film flows down an inclined ramp. We provide a set of equations, namely Whitham's modulation equations, that determines the behaviour of low frequency perturbations of periodic wave trains. As a byproduct, we relate the spectral stability in the small wavenumber regime to properties of the modulation equations. This stability is always critical since $0$ is a $0$-Floquet number eigenvalue associated to translational invariance.

\section{Introduction}

We study perturbations of periodic travelling wave solutions to the following generalized Kuramoto-Sivashinski (gKS) equations 
\begin{equation}\label{gKSi}
\displaystyle
\partial_t u+6u\partial_x u+\partial_{x}^3u+\delta\Big(R\partial_x^2 u+\partial_x^4u\Big)\ =\ 0
\end{equation}
where $u(x,t)\in\R$ is a heigth at place $x\in\R$ and time $t>0$. Actually this equation describes approximately long and small surface waves of two dimensional incompressible viscous fluid down an inclined plane. Here $\delta R$ ($\delta>0$ and $R>0$) measures the deviation of the Reynolds number from the critical Reynolds number above which long wave perturbations are spectrally unstable. Since $R>0$, steady flows are unstable and periodic travelling waves, so called roll-waves, appear. 

This transition to instability is proved to occur also in shallow water equations (SW): 
\begin{equation}\label{sv}
\left\{
\begin{array}{rcl}
\displaystyle
\partial_t h\ +\ \partial_x (hu)&=&0\\
\displaystyle
\partial_t (hu)\ +\ \partial_x\left(hu^2+\frac{h^2}{2F^2}\right)&=&h-u^2+R_e^{-1}\partial_x(h\partial_x u)
\end{array}
\right.
\end{equation}
(here $h(t,x)\in\R$ is an heigth, $u(t,x)\in\R$ a velocity, $F>0$ a Froude number and $R_e>0$ a Reynolds number).
If $F>2$, the stationnary solutions are unstable and roll-waves also appear. These roll-waves are observed both experimentally and numerically.

Both in \eqref{gKSi} and \eqref{sv}, roll-waves are proved to exist through a Hopf bifurcation argument and a whole family of periodic solutions exists, which ends up with a solitary wave. This family is parametrized by two quantities, period and either spatial mean heigth for \eqref{gKSi} or relative discharge rate for \eqref{sv}. Under some conditions, a similar situation is also proved to exist in the Navier-Stokes equations (NS) for an incompressible fluid down a ramp: indeed, taking into account capillarity and under a spectral assumption, namely a Hopf bifurcation scenario, it is proved that small amplitude periodic travelling wave solutions to (NS) exist when $R_e>5\,{\rm cotan}(\theta)/6$, $\theta$ being the angle of the slope \cite{Ni_Te_Yo}. 

Actually the (gKS) equation is a generic equation which describes this transition to instability in the regime of small amplitude long waves. One can derive formally (gKS) either from (SW) or (NS): see e.g. \cite{Win} for a derivation from (NS) and \cite{Yu_Yang} for a derivation from (SW). These formal derivations are strongly connected to the issue of the \emph{rigorous} derivation of the shallow water equations from the Navier Stokes equations in the case of viscous fluids flowing down a ramp. Assuming capillarity, the (SW) equations are derived rigorously from (NS) whenever steady solutions are stable \cite{Br_No}. The (gKS) equations is the reference equation to describe small amplitude solutions of (NS) or (SW) when the flow is weakly unstable. In the stable regime, both the shallow water equations and (NS) can be approximately reduced to a single viscous Burgers equation (or Benney equation): see \cite{Uecker} for the reduction from (NS), as for (SW) the limit is a standard relaxation limit in balance laws. Somehow this reduction already provides us with a validation of (SW) in the stable regime. Yet solutions to the Burgers equation blow up when $R_e$ exceeds a critical value and (gKS) is then needed.

Stability of roll-waves and description of the dynamic around such wave trains are interesting intricated questions. Some information may already be deduced from the fact that roll-waves emerge from unstable steady solutions. At the onset of roll-waves, that is in the small amplitude, the wave train is close to a stationnary solution and so is expected to be spectrally unstable. For the same reason, the solitary wave is also expected to be spectrally unstable, so that only a range of periodic roll-waves could be stable. This would explain the coalescence dynamic that is usually observed (see e.g. \cite{Bar}). To go beyond this general heuristic discussion, we summarize here some of the results obtained for the (SW) equations in \cite{Ba_Jo_No_Ro_Zu_2,Ba_Jo_No_Ro_Zu_1,Ba_Jo_Ro_Zu}. It is always a non trivial task to carry out the spectral analysis of equations linearized about periodic travelling waves. Yet, using perturbation methods or reducing to numerical computations in a finite box through energy estimates or stable-unstable tracking, solitary waves and small amplitude wave trains are there proved to be spectrally unstable and a range of periodic roll-waves is indeed proved to be stable. The low frequency instability is a strong instability for small amplitude (the spectrum crosses transversally the imaginary axis), whereas, for solitary waves, essential spectrum is unstable (corresponding to steady states spectrum) and point spectrum is weakly unstable (the spectrum is tangent to the imaginary axis but on the unstable side). Note that the spectral stability of some periodic wave trains may be used to prove also their nonlinear asymptotic stability \cite{Jo_Zu_No}.

Complementary to direct numerical computations, another approach to study the stability of periodic wave trains is to carry out a long wavelength analysis in the neighbourhood of roll-waves since this is the particular regime governing the return to equilibrium. Slow modulation of roll-waves are then proved to be (still) described locally by roll-waves whose parameters evolve on a slow scale according to a set of first order partial differential equations, so called Whitham's averaged equations. This approach involves at least two steps, the derivation of averaged equations and the proof of a connection between spectral proprerties of both the modulation system and the original equations. Obviously the spectral step yields then some necessary conditions,  expressed through averaged quantities, for the stability of periodic travelling waves. Such an approach was proposed for conservation laws in \cite{Serre} where it is proved that the leading term, in the expansion about the origin, of the corresponding Evans function is given by the dispersion relation associated to linearized Whitham's equations. For (SW), an Evans function approach to stability was also introduced in \cite{N1} and extensively used in \cite{Ba_Jo_No_Ro_Zu_2,Ba_Jo_No_Ro_Zu_1,Ba_Jo_Ro_Zu}, and his relation with the dispersion of the Whitham's system also established in \cite{No_Ro}. Yet, in \cite{No_Ro} a Bloch transform approach is shown to give a more natural and powerful way of achieving the spectral step. Actually, in \cite{No_Ro}, the authors obtained, as in \cite{Serre} for conservation laws, that the tangency to the imaginary axis of low frequency spectral curves is descibed by the hyperbolicity of a first order modulation system but also that the curvatures of this curves are related to the parabolicity of a second order modulation system also derived in \cite{No_Ro}. Moreover, the Bloch transform approach gives almost readily a relation not only between eigenvalues but also eigenvectors. This is precisely the keystep in order to perform a third step in the modulation analysis : to validate at a nonlinear level the modulation equations as providing a good approximation of low-frequency perturbations of wave trains as is explained in \cite{D3S} for reaction-diffusion systems (in cases where modulation yields a single scalar equation and not a system). In \cite{No_Ro} the authors proved such a nonlinear validation of the first order modulation system when its hyperbolicity is met. An expansion of eigenvectors would certainly also be useful for a justification of second order modulation systems and constructions of generalized shocks solutions corresponding to modulated roll-waves with shocks in parameters as performed for scalar modulation equations in reaction-diffusion systems \cite{D3S}.

In this paper, we carry out the derivation and spectral steps of a similar analysis for (gKS). For fixed $\delta>0$, the set of periodic traveling waves is two dimensional (up to translations) and it is proved to exist through a Hopf bifurcation argument. This set is parametrized by the wavenumber and spatial mean. In the low frequency regime, we compute a set of two equations that governs slow modulation of the local wavenumber and local spatial mean. Similarly to roll-waves in (SW), we obtain "inviscid" (first order) and "viscous" (second-order) models and relate these systems to the spectral stability of periodic wave trains just as in \cite{No_Ro}.

However, recall form the derivation of (gKS) from (NS) and (SW) that the parameter $\delta$ should be small so that the domain of validity of Whitham's equations, whose size is measured as the size of the allowed perturbation characteristic wavenumber in the modulation, shrinks to $0$. One has to carry out a new modulation analysis with characteristic wavenumbers $\varepsilon$ comparable to $\delta$. We do this setting $\delta=\bar\delta \varepsilon$ and we show in this case that the Whitham set of equations is composed of three first order partial differential equations with a balance term which takes into account the dissipation of (gKS). The number of equations, three, is easily understood since, in the limit $\delta\to 0$, (gKS) is a perturbation of the Korteweg-de Vries (KdV) equation which possesses a three dimensional set of periodic wave trains. We will connect this set of modulation equations with the spectral stability of periodic wave trains in the regime of wavenumber of order $\delta$. We also relate the modulation system with the one for (KdV) (in the limt $\bar\delta\to0$) and the one obtained for a fixed $\delta$ (in the limit $\bar\delta\to\infty$).

The paper is organized as follows. In section \ref{sec2}, we recall the structure of periodic wavetrains solutions to (gKS) both in the case $\delta>0$ fixed and in the limit $\delta\to 0$. In section \ref{sec3}, we consider the case $\delta>0$ fixed, we first derive "inviscid" and "viscous" modulation equations, then carry out a direct spectral analysis of (gKS) and relate it with a spectral analysis of the modulation systems. In section \ref{sec4}, we consider the case $\delta\to 0$, compute the set of Whitham's equations, carry out a direct spectral analysis in the regime of frequencies of order $\delta\to0$ and again show the relation with the Whitham's systems.

\section{\label{sec2} Periodic travelling waves in (gKS)}

Up to changes
$$
t=R^{3/2}t,\ x=R^{1/2}x,\ u=Ru, \delta=R^{1/2}\delta
$$
one may assume\footnote{One could also eliminate $6$ by rescaling $u$ again but we keep it to fit with equations as written in \cite{Bar}.} $R=1$ and consider (gKS) equation in the form
\begin{equation}\label{gKS}
\displaystyle
\partial_t u+6 u\partial_x u+\partial_{x}^3 u+\delta\big(\partial_x^2 u+\partial_x^4 u\big)=0\ .
\end{equation}
For the moment we search for periodic travelling waves, with profile $U$ a periodic function, without fixing the period of the profile, and thus for $u$ in the form $u(x,t)=U(x-ct)$. This yields
\begin{equation}\label{ep}
\displaystyle
-cU'+6UU'+U''+\delta\big(U''+U''''\big)\ =\ 0\ .
\end{equation}
Integrating once, we search for a periodic profile $U$, a speed $c$ and a constant $\bq$ such that  
\begin{equation}\label{epI}
\displaystyle
-cU+3U^2+U''+\delta\big(U'+U'''\big)\ =\ \bq\ .
\end{equation}

For arbitrary $\delta>0$, the existence of periodic solutions to \eqref{ep} may be obtained through a Hopf bifurcation argument. Equation (\ref{epI}) possesses two stationnay solutions $U_{-}(c,\bq)<U_{+}(c,\bq)$ such that $\displaystyle 3U_{\pm}^2-cU_{\pm}=\bq$. Linearizing (\ref{epI}) at $U=U_{\pm}$ yields  
\begin{equation}\label{ep_lin}
\displaystyle
(6U_{\pm}-c)\tilde{U}+\tilde{U}''+\delta\big(\tilde{U}'+\tilde{U}'''\big)\ =\ 0\ .
\end{equation}
The assocatied characteristic equation is $\displaystyle\delta\big(\lambda^3+\lambda\big)+\lambda^2+(6U_{\pm}-c)=0$. It is easily seen that the stationnary solution undergoes a Hopf bifurcation whenever $6U_{\pm}-c=1$. This can only happens for $U=U_+$ and the equation settles a $c_{Hopf}(\bq)$. Denote $V=U-U_+$ and $6U_+-c=1+\mu$ with $\mu\ll 1$. Equation \eqref{epI} then reads
$$
\displaystyle
V(3V+1+\mu)+V''+\delta(V'+V''')\ =\ 0\ .
$$
For $\mu$ sufficiently small on one side of $0$, there is a family $V(\mu)$ of amplitude $\mathcal{O}(\sqrt{|\mu|})$ and frequency $k(\mu)$, with $k(0)=1$ and $k'(0)\neq 0$. 

As a result, one obtains a two dimensional manifold of periodic traveling waves (identified when coinciding up to translation) parametrized by the wavenumber $k$ and $\bq$. Note that in the small amplitude regime one may also parametrize by $k$ and the spatial mean of the solution. Indeed, there stands $U=<\!U\!>+\mathcal{O}(\sqrt{|\mu|})$, thus since $\displaystyle \bq=3<\!U^2\!>-c<\!U\!>$, one finds $\displaystyle 3<\!U\!>^2-c<\!U\!>=\bq+\mathcal{O}(|\mu|)$. At last recall that, for $\mu$ small enough, $6<\!U\!>-c\neq 0$ so that one can switch from $\bq$ to $M=<\!U\!>$ in the neighbourhood of the Hopf bifurcation. As a result one may prove the following proposition.
\begin{proposition}
Let $M\in\R$ be fixed and $k<1$ be such that $1-k$ is small. Then there exist a unique $c(k,M)$ and a unique $\bq(k,M)$ such that there exist a $1$-periodic $U(\,\cdot\,;k,M)$ solution to 
\begin{equation}
\label{epIk}
\displaystyle
k\left(3U^2-cU\right)+k^3U''+\delta\left(k^2U'+k^4U'''\right)\ =\ \bq,\quad <U>=M\ .
\end{equation}
Moreover this solution is unique up to translation.
\end{proposition}

Maybe the last thing we should say about the proof of this proposition is why the bifurcation occurs for $k<1$ (rather than $k>1$). A quick way to see this is to multiply the first equation in \eqref{epIk} by $U'$ and integrate over a period. This yields $<(U')^2>=k^2<(U'')^2>$ thus $k\leq1$. Note also that above and through the text we use $<f>$ to denote the average of $f$ over one period ; for $1$-periodic functions this is just $\int_0^1f$.

Let us now explain how to reduce the search of solutions to \eqref{epIk} to the case $M=0$. Inserting the \emph{ansatz} 
\begin{equation}\label{epM}
U=M+U^{(0)}
\end{equation} into \eqref{epIk} yields $<U^{(0)}>=0$ and again
\begin{equation*}
\displaystyle
-c^{(0)}U^{(0)}+3U^{(0)}{}^2+k^2U^{(0)}{}''+\delta\left(k U^{(0)}{}'+k^3U^{(0)}{}'''\right)\ =\ \bq^{(0)}
\end{equation*}
with
\begin{equation}\label{epMparam}
c^{(0)}\ =\ c-6M\,,\quad\bq^{(0)}\ =\ \bq+cM-3M^2\,.
\end{equation}
Note that in the $(k,M)$ parametrization of profiles, none of the above quantities $c^{(0)}$, $U^{(0)}$, $\bq^{(0)}$ depends on $M$. From now on all quantities with a $(0)$ in superscript will refer to the zero mean problem.

For applications, we usually assume that $\delta$ is small. Then equation \eqref{gKS} is a singular perturbation of (KdV) equation whose periodic solutions are described with the help of elliptic functions. In the singular limit $\delta\to 0$, periodic solutions to \eqref{ep} were proved in \cite{Er_Mc_Ro} to be close to periodic solutions to (KdV) and an expansion with respect to $\delta$ was obtained. Moreover, in \cite{Bar}, a formal spectral analysis was carried out, but only in the context of perturbations with zero mean. This does induce a lack of generality. Here, we will provide a spectral stability analysis for \emph{arbitrary} (bounded) perturbations.

We provide now expansions of roll-waves profiles in the $\delta\to0$ limit. As $\delta\to 0$, one may expand $c^{(0)}$, $U^{(0)}$, $\bq^{(0)}$ as
$$
\displaystyle
c^{(0)}=\tilde{c}_0+\delta\tilde{c}_1+\mathcal{O}(\delta^2),\quad \bq^{(0)}=\tilde{q}_0+\delta\tilde{q}_1+\mathcal{O}(\delta^2),\quad U^{(0)}=\tilde{U}_0+\delta\tilde{U}_1+\mathcal{O}(\delta^2).
$$
The $1$-periodic solution to \eqref{epI} for $\delta=0$  and $<\tilde{U}_0>=0$ is given by
\begin{equation}\label{per_kdv}
\begin{array}{rcl}
\displaystyle
\tilde{U}_0(\xi)&=&8K(p)^2\,k^2\,\left({\rm dn}^2(2k\,K(p)\xi)-\frac{E(p)}{K(p)}\right),\\
\displaystyle
\tilde{c}_0&=&16K(p)^2{\pi}\,k^2\left(2-p^2-3\frac{E(p)}{K(p)}\right),\\
\displaystyle
\tilde{q}_0&=&4\left(2K(p)\,k\right)^4\left(-3\frac{E(p)^2}{K(p)^2}+2(2-p^2)\frac{E(p)}{K(p)}+p^2-1\right),
\end{array}
\end{equation}
with $p\in]0, 1[$, $K(p)$ and $E(p)$ elliptic integrals of first kind and second kind and $\textrm{dn}(\,\cdot\,)$ Jacobi's delta function with modulus $p$. Up to this order, there is no selection of a particular wave train ($p$ is arbitrary) and the manifold of periodic solutions (identified when coinciding up to translation) would be of dimension $2$, whereas it was $1$-dimensional for $\delta>0$ fixed when restricting to zero mean solutions. To recover the same dimension, one has to compute an expansion of the solution to the next order: one finds
\begin{equation}\label{dl1}
\displaystyle
k^2\tilde{U}_1'''+(6\tilde{U}_0\tilde{U}_1-\tilde{c}_0\tilde{U}_1)'-\tilde{c}_1\tilde{U}_0'+k\tilde{U}_0''+k^3\tilde{U}_0''''\ =\ 0\ .
\end{equation}
The linear operator $\displaystyle l=k^2\frac{d^3}{dx^3}+\frac{d}{dx}\left((6\tilde{u}_0-\tilde{c}_0)\,\cdot\,\right)$ is Fredholm of index $0$ and $1$ and $\tilde{U}_0$ span the kernel of its adjoint so that one can readily deduce that equation \eqref{dl1} has a solution provided that the following compatibility condition is satisfied, $\displaystyle <(\tilde{U}'_0)^2>=k^2<(\tilde{U}''_0)^2>$. This equation yields the selection criterion
\begin{equation}\label{eq9}
\frac{1}{k^2}\ =\ 16\mathcal{F}(p)
\end{equation}
where
\begin{equation}
\begin{array}{rcl}
\displaystyle
\mathcal{F}(p)&=&3E(p)K(p)+(p^2-2)K(p)^2-\dfrac{6N(p)}{7D(p)}\ ,\\[1ex]
\displaystyle
N(p)&=&7(1-p^2+p^4)E(p)^2-(10-15p^2+13p^4-4p^6)E(p)K(p)\\[1ex]
\displaystyle
&&+(3-6p^2+5p^4-2p^6)K(p)^2\ ,\\[1ex]
\displaystyle
D(p)&=&2(1-p^2+p^4)E(p)/K(p)-(2-3p^2+p^4)\ .
\end{array}
\end{equation}
In order to determine $\tilde{c}_1$, one has to consider higher order corrections to $U^{(0)}$: in fact,  $\tilde{c}_1$ is determined through a solvability condition on the equation for $\tilde{U}_2$. This yields $\tilde{c}_1=0$ (see \cite{Er_Mc_Ro} for more details).
As a consequence, coming back to the more general case where $M$ may not vanish, we have also obtained a two dimensional manifold of (asymptotic) periodic solutions (identified when coinciding up to translation) para\-metrized by spatial mean $M$ and wavenumber $k$ (or alternatively the parameter $p\in]0, 1[$). Note that the limit $k\to 0$ (i.e. $p\to 1$) corresponds to a soliton and $k\to 1$ (i.e. $p\to 0$) corresponds to small amplitude solutions (or equivalently  to the onset of the Hopf bifurcation branch).

In what follows, we establish Whitham's equations that describes modulation of these wave trains through slow evolution of their parameters. For fixed $\delta>0$, this will give us a system of two partial differential equations. To analyze the limit $\delta\to 0$ and obtain a consistant Whitham's modulation theory, one has to do something slightly different and choose a characteristic perturbation wavenumber $\varepsilon$ of order $\mathcal{O}(\delta)$. The equation (gKS) is then a perturbation of (KdV) and we obtain a set of three partial differential equations which are balance laws but with a source term. As $\delta/\varepsilon\to\infty$, this source term is stiff and the system relaxes to a system of two partial differential equations similar to the one written for fixed $\delta>0$. 

Let us also introduce the notation $\Omega=-kc$ for the rest of this work.

\mathversion{bold}
\section{\label{sec3} Stability of wave trains for fixed $\delta>0$}
\mathversion{normal}

\subsection{Modulation equations}
\subsubsection{First order Whitham's equations}

Let us denote $\varepsilon$ the characteristic wavenumber of perturbations. To study low frequency perturbations on scales $(x/\varepsilon,t/\varepsilon)$ and obtain an averaged modulated behaviour, we first rescal variables as $(X,T)=(\varepsilon x,\varepsilon t)$. Then equation \eqref{epI} turns into
\begin{equation}\label{gKS_m}
\displaystyle
\partial_{T} u+6u\partial_X u+\varepsilon^2\partial_X^3 u+\delta\big(\varepsilon\partial_X^2 u+\varepsilon^3\partial_X^4 u)\ =\ 0\ .
\end{equation}
In what follows, we will assume $\delta>0$ fixed (but small). We search for solutions to \eqref{gKS_m} in the (formal) form
$$
\displaystyle
u(X,T)=\sum_{k=0}^{\infty}\varepsilon^k u_k\left(\frac{\varphi(X,T)}{\varepsilon};X,T\right)\quad\textrm{with}\quad\varphi(X,T)=\sum_{j=0}^{\infty}\varepsilon^j\varphi_j(X,T)\ ,
$$
where $u_i(y;X,T)$ are $1$ periodic functions in $y$. Inserting this ansatz into \eqref{gKS_m} and setting $\Omega_0=\partial_{T}\phi_0$ and $k_0=\partial_X\phi_0$, one finds at the $\mathcal{O}(\varepsilon^{-1})$ order
\begin{equation}\label{me_0}
\displaystyle
\Omega_0\partial_y u_0+6k_0u_0\partial_y u_0+k_0^3\partial_y^3 u_0+\delta\left(k_0^2\partial_y^2 u_0+k_0^4\partial_y^4 u_0\right)\ =\ 0\ .
 \end{equation}
Denoting $M_0=<u_0>$, equation \eqref{me_0} implies $\Omega_0=\Omega(k_0,M_0)$ and (up to some translation) 
$$
u_0(y;X,T)\ =\ U(k_0(X,T),M_0(X,T))(y)\ .
$$
Compatibility condition $\displaystyle \partial_T\partial_X\varphi_0=\partial_X\partial_T\varphi_0$ yields the first equation of a Whitham's system
\begin{equation}\label{w1_0}
\displaystyle
\partial_T k_0-\partial_X\left(\Omega(k_0,M_0)\right)\ =\ 0\ .
\end{equation}
Recalling $\Omega=-kc$, this equation is also written as
\begin{equation}\label{w1_00}
\displaystyle
\partial_Tk_0+\partial_X\left(k_0c(k_0,M_0)\right)\ =\ 0\ ,
\end{equation}
Now we identify $\mathcal{O}(\varepsilon^{0})$ terms. Setting $\Omega_1=\partial_T\varphi_1$ and $k_1=\partial_X\varphi_1$, one finds
\begin{equation}\label{me+0}
{\setlength\arraycolsep{1pt}
\begin{array}{rcl}
\displaystyle\!\!\!
\LL_0\,u_1&+&\partial_Tu_0+\partial_X(3u_0^2)+3\left(k_0^2\partial^3_{Xyy}u_0+k_0\partial_Xk_0\partial_y^2 u_0\right)\\[1ex]
\displaystyle
&+&\delta\left(2k_0\partial_{Xy}u_0+\partial_{X}k_0\partial_y u_0+4k_0^3\partial^4_{Xyyy}u_0+6k_0^2\partial_X k_0\partial_{y}^3u_0\right)\\[1ex]
\displaystyle
&+&\Omega_1\partial_y u_0+k_1\left(\partial_y(3u_0^2)+3k_0^2\partial^3_y u_0+\delta(2k_0\partial_y^2 u_0+4k_0^3\partial_y^4 u_0) \right)=0
\end{array}}
\end{equation}
where $\LL_0$ denotes the linear differential operator defined by
\begin{equation*}
\begin{array}{l}
[\LL_0 f](y;X,T)\ =\ \left[\LL_{k_0(X,T),M_0(X,T)}f(\,\cdot\,;X,T)\right](y)\ ,\\[1ex]
\displaystyle
\LL_{k,M}f=k\frac{d}{dy}\left(\left((6U(k,M)-c(k,M)\right)f\right)+k^3\frac{d^3f}{dy^3}+\delta\left(k^2\frac{d^2f}{dy^2}+k^4\frac{d^4f}{dy^4}\right)\ ,
\end{array}
\end{equation*}
$\LL$ being thus the operator associated to the linearization about $U(k,M)$ of the profile equation in \eqref{epIk}  differentiated once. The kernel of $\LL=\LL_{k,M}$ is one dimensional and spanned by $U'(k,M)$. Indeed, $0$ is a semi-simple eigenvalue of $\LL$ of multiplicity $2$ and $\partial_M U(k,M)=1$ is a generalized eigenvector since $\LL[1]=6kU'$. The adjoint of $\LL$, $\LL^{ad}$, is Fredholm of index $0$ with a one dimensional kernel\footnote{The fact that $1$ lies in the kernel is directly related to the fact that (gKS) is a conservation law.} spanned by $1$. As a result, equation \eqref{me+0} is solvable if 
\begin{equation}\label{w2_00}
\displaystyle
\partial_T<u_0>+\partial_X<3u_0^2>\ =\ 0.
\end{equation}
\noindent
Note that using the notations \eqref{epMparam} corresponding to the decomposition \eqref{epM}, one has 
$$
\displaystyle
c\ =\ 6M_0+c^{(0)}(k_0)\quad\textrm{and}\quad<3u_0^2>\ =\ \bq+cM_0\ =\ \bq^{(0)}(k_0)+3M_0^2\ .
$$ 
This yields the following Whitham's system
\begin{equation}\label{w12_0}
\displaystyle
\partial_T k_0+\partial_X\left(k_0\, c(k_0,M_0)\right)=0,\quad \partial_T M_0+\partial_X\left(\bq^{(0)}(k_0)+3M_0^2\right)=0.
\end{equation}
 
This is the classical Whitham's system in the case of (gKS). Later on, we will relate the hyperbolicity of this system with a stability index that is used to give necessary conditions for the spectral stability of periodic wave trains.

\subsubsection{Second order Whitham's equations}

We now proceed with the modulation expansion and find an evolution system for $k_1$ and $M_1=<u_1>$. This system will be useful to complete the spectral stability analysis of periodic wave trains. We first rewrite equation \eqref{me+0} in a more convenient way. Differentiating \eqref{epIk} with respect to $k$ and $M$ yields
\begin{equation}\label{dku}
\begin{array}{l}
\displaystyle
\LL[\partial_k U]=-\partial_k\Omega\,U'-\left((3U^2)'+3k^2U'''+\delta(2kU''+4k^3U'''')\right)\,,\\[1ex]
\displaystyle
\LL[\partial_M U]=-\partial_M\Omega\,U'\,,\quad<\partial_k U>=0\,,\quad<\partial_M U>=1\,.
\end{array}
\end{equation}
As a consequence, denoting by a subscript $0$ the evaluation in $(k_0,M_0)$ (so that for instance $u_0=U_0$), equation \eqref{me+0} may be turned into
\begin{equation}\label{me+01}
\displaystyle
\LL_0[u_1-k_1\partial_k U_0]+\left(\Omega_1-k_1\partial_k\Omega\right)\partial_y u_0+\partial_T u_0+\partial_X(3u_0^2)+R_0
\end{equation}
with $R_0$ defined as 
\begin{equation}
\begin{array}{rcl}
\displaystyle
R_0&=&3\left(k_0^2\partial^3_{Xyy}u_0+k_0\partial_Xk_0\partial_y^2 u_0\right)+\delta\left(2k_0\partial_{Xy}u_0+\partial_{X}k_0\partial_y u_0\right)\\[1ex]
\displaystyle
&&+\delta\left(4k_0^3\partial^4_{Xyyy}u_0+6k_0^2\partial_X k_0\partial_{y}^3u_0\right)\ .
\end{array}
\end{equation}
Note that $R_0$ depends only on $(k_0,M_0)$ and has zero mean. From \eqref{w2_00}, we know that there is a unique $f_0$ satisfying $\displaystyle \LL_0[f_0]=\partial_Tu_0+\partial_X(3u_0^2)+R_0$ and $<f_0,U_0'>=0$. Equation \eqref{me+01} now reads
\begin{equation}\label{me+02}
\displaystyle
\LL_0[u_1-k_1\partial_k U_0+f_0-<f_0>]+\left(\Omega_1-k_1\partial_k\Omega+6k_0<f_0>\right)U'_0\ =\ 0\ .
\end{equation}  
Let us denote by $M_1$ the function of $(X,T)$ satisfying 
\begin{equation}\label{Om1}
\displaystyle \Omega_1-k_1\partial_k\Omega_0-M_1\partial_M\Omega_0+6k_0<f_0>\ =\ 0\ . 
\end{equation}
Then equation \eqref{me+02} is equivalent to  
\begin{equation}\label{u1}
\displaystyle u_1=M_1\partial_M U_0+k_1\partial_k U_0+<f_0>-f_0+A_1 U_0'
\end{equation}
for some $A_1$ (independent of $y$). Note that $M_1$ appears to be $M_1=<u_1>$. We now write a set of evoultion equations for $(k_1,M_1)$. We first obtain an equation that governs the time evolution of $k_1$ by using the compatibility condition $\partial_T\partial_X\varphi_1=\partial_X\partial_T\varphi_1$ and  equation \eqref{Om1},
\begin{equation}\label{w1_1}
\displaystyle
\partial_T k_1-\partial_X\Big(k_1\frac{\partial\Omega}{\partial k}+M_1\frac{\partial\Omega}{\partial M}\Big)
\ =\ -\,\partial_X<6k_0\,f_0>.
\end{equation}
\noindent
In order to obtain an equation that governs the evolution of $M_1$, we have to consider $\mathcal{O}(\varepsilon)$ terms in \eqref{gKS_m},
\begin{equation}\label{me_1}
\displaystyle
\partial_T u_1+\partial_X\big(6u_0\,u_1\big)+\delta\partial_{XX}u_0\ =\ \partial_y\left(\cdots\right)\ .
\end{equation}
The exact form of the right hand side in \eqref{me_1} is not needed here to compute the set of equations for $(k_1,M_1)$. Indeed, equation \eqref{me_1} imposes the compatibility condition
\begin{equation}\label{w1_20}
\displaystyle
\partial_T<u_1>+\partial_X<6u_0\,u_1>+\delta\partial_{XX}<u_0>\ =\ 0\ .
\end{equation}
Inserting \eqref{u1} into \eqref{w1_20} yields
\begin{equation}\label{w1_2}
\begin{array}{rcl}
\displaystyle
\partial_T M_1&+&\partial_X\left(\textrm{d}\left[<3U^2>\right][k_1, M_1]\right)\\[1ex]
&=&\partial_X<6u_0(f_0-<f_0>)>-\delta\partial_{XX}M_0
\end{array}
\end{equation}
where $\textrm{d}$ yields differential with respect to parameters $(k,M)$. This almost completes the derivation of a viscous Whitham's system. Indeed, let us introduce $k=k_0+\varepsilon k_1$ and $M=M_0+\varepsilon M_1$. Next, consider \eqref{w1_00}$+\varepsilon$\eqref{w1_1} on the one hand and \eqref{w2_00}$+\varepsilon$\eqref{w1_2}. Then, neglecting $\mathcal{O}(\varepsilon^2)$ terms, these equations are written as a system of viscous conservation laws with $\mathcal{O}(\varepsilon)$ viscous terms. We shall relate these viscous equations with the curvature of the spectral curves associated to the stability of periodic wave trains in the small wavenumber regime.

\subsection{Spectral stability of periodic wave trains}

\subsubsection{Bloch analysis of spectral stability}

In this section, we analyse the spectral stability of a periodic wave train $U(\,\cdot\,;\ks,\Ms)$ Denoting by a star subscript the evaluation in $(\ks,\Ms)$, in an adapted co-moving frame equation \eqref{gKS} linearized about $\Us$ reads
\begin{equation}\label{gKS_l}
\displaystyle
\partial_T u+\ks\partial_x\left((6\Us-\cs)u\right)+\ks^3\partial_x^3u+\delta\left(\ks^2\partial_x^2{u}+\ks^4\partial_{x}^4u\right)=0.
\end{equation}
This is shortly written as $\partial_t u+\LLs u$ where $\LLs=\LL_{\ks,\Ms}$. We carry out a Bloch analysis of this problem and search solutions to \eqref{gKS_l} in form $u(x,t)=e^{\lambda t+\nu x}\hat{u}_{\nu}(x)$, with $\hat{u}_{\nu}$ a $1$-periodic function and $\nu\in i\R$. For writting convenience, we will drop the $\nu$ subscript indicating $\nu$ dependence. As announced we consider
\begin{equation}\label{gKS_ls}
\begin{array}{rcl}
\displaystyle
\ks\left(\frac{d}{dx}+\nu\right)\left((6\Us-\cs)u\right)
&+&\displaystyle
\ks^3\left(\frac{d}{dx}+\nu\right)^3\hat{u}\\[1ex]
\displaystyle
+\ \delta\Big[\ks^2\left(\frac{d}{dx}+\nu\right)^2\hat{u}
&+&\displaystyle
\ks^4\left(\frac{d}{dx}+\nu\right)^4\hat{u}\Big]
\ =\ -\lambda\,\hat{u}\ .
\end{array}
\end{equation}
The value $\lambda=0$ is an eigenvalue of $\LLs$ associated to the translation mode $\hat{u}=U'$ thus corresponding to $\nu=0$. We compute the spectrum in the regime of large wavelength perturbations $|\nu|\ll 1$ and expand solutions to equation \eqref{gKS_ls} as
$$
\displaystyle
\lambda\ =\ \nu\lambda_0+\nu^2\lambda_1+\mathcal{O}(\nu^3),
\qquad
\hat{u}\ =\ \hat{u}_0+\nu\hat{u}_1+\nu^2\hat{u}_2+\mathcal{O}(\nu^3).
$$
We first identify $\mathcal{O}(\nu^0)$ terms in equation \eqref{gKS_ls} and find $\LLs\hat{u}_0=0$ so that  $\hat{u}_0=k_0\Us'/\ks$ for some constant $k_0$. Note that $\hat{u}_0$ is mean free.

Next, we identify $\mathcal{O}(\nu^1)$ terms and obtain  
\begin{equation*}
\displaystyle
\LLs\hat{u}_1+\frac{\lambda_0k_0}{\ks}\Us'+k_0\left((6\Us-\cs)\Us'+3\ks^2\Us'''+\delta\left(2\ks\Us''+4\ks^3\Us''''\right)\right)=0
\end{equation*}
thus, using \eqref{dku} and denoting $M_0$ the constant satisfying
\begin{equation}\label{cond_10}
\displaystyle
\frac{\lambda_0\,k_0}{\ks}+\ks\partial_k\cs\,k_0+\ks\partial_M\cs\,M_0\ =\ 0\ ,
\end{equation}
we are left with
\begin{equation}\label{eq_u1}
\displaystyle
\LLs[\hat{u}_1-k_0\partial_k \Us-M_0\partial_M\Us]\ =\ 0\ .
\end{equation}
But equation \eqref{eq_u1} is equivalent to
\begin{equation}\label{bar_u1}
\displaystyle
\hat{u}_1\ =\ k_0\,\partial_k \Us+M_0\,\partial_M\Us+\frac{k_1}{\ks}\Us'
\end{equation} 
for some constant $k_1$. Note that $M_0=<\hat{u}_1>$. In order to obtain a full set of equations for $(k_0,M_0)$, one has to consider $\mathcal{O}(\nu^2)$ terms. We find that $\hat{u}_2$ should satisfy
{\setlength\arraycolsep{1pt}
\begin{eqnarray}
\displaystyle
\LLs[\hat{u}_2]&+&\lambda_0\hat{u}_1+\ks\left((6\Us-\cs)\hat{u}_1+3\ks^2\frac{d^2\hat{u}_1}{dx^2}+\delta\left(2\ks\frac{d\hat{u}_1}{dx}+4\ks^3\frac{d^3\hat{u}_1}{dx^3}\right)\right)\nonumber\\
\displaystyle
\label{eq_u2}
&+&\lambda_1\hat{u}_0+3\ks^3\frac{d\hat{u}_0}{dx}+\delta\left(\ks^2\hat{u}_0+6\ks^4\frac{d^2\hat{u}_0}{dx^2}\right)
\ =\ 0\ .
\end{eqnarray}}
Equation \eqref{eq_u2} has a solution if and only if 
$$
\displaystyle
\lambda_0<\hat{u}_1>+\ks<\big(6\Us-\cs\big)\hat{u}_1>\ =\ 0\ .
$$
Substituting \eqref{bar_u1} into this equation yields the equation
\begin{equation}\label{cond_20}
\displaystyle
\ks\partial_k[<3U^2>]_\star\,k_0+\left(\lambda_0-\ks\cs+\ks\partial_M[<3U^2>]_\star\right)\,M_0
\ =\ 0\ .
\end{equation}
Setting $\lambda_0=\ks\tilde{\lambda}_0$, the linear system (\ref{cond_10},\ref{cond_20}) has a solution if and only if
\begin{equation}\label{disp_0}
\displaystyle
\left|\begin{array}{cc} \displaystyle \tilde{\lambda}_0+\ks\partial_k \cs &\displaystyle  \ks\partial_M\cs\\
\displaystyle \partial_k[<3U^2>]_\star&\quad \displaystyle \tilde{\lambda}_0-\cs+\partial_M[<3U^2>]_\star\end{array}\right|
\ =\ 0\ .
\end{equation}
As a result, one obtains a necessary condition for the stability of periodic wave trains. Indeed, recall that $\lambda$ expands as $\lambda=\nu\ks\tilde{\lambda}_0+\mathcal{O}(\nu^2)$ and $\nu\in i\R$ so that for the periodic wave train to be stable, it is necessary that $\tilde{\lambda}_0\in\mathbb{R}$. In next section, we will find out this former condition also expresses as hyperbolicity of the first order Whitham's system.

Again we proceed with the computation of the expansion of eigenvalues.
Let us first introduce $\hat{v}^k$ and $\hat{v}^M$ defined as
\begin{equation*}
\begin{array}{rcl}
\displaystyle
\hat{v}^k&=&(\tilde{\lambda}_0-\cs)\,\partial_k\Us+3\,\partial_k\left(U^2-<U^2>\right)_\star
+3\,\ks(\Us+\ks\partial_k\Us)''\\[1ex]
\displaystyle
&&+\delta\left(\Us'+6\,\ks^2\Us'''+2\,\ks(\partial_k\Us)'+4\,\ks^3(\partial_k\Us)'''\right)\ ,\\[1ex]
\displaystyle
\hat{v}^M&=&(\tilde{\lambda}_0-\cs)\,\partial_M\Us+3\,\partial_M\left(U^2-<U^2>\right)_\star\ .
\end{array}
\end{equation*}
Notice that $\hat{v}^k$ and $\hat{v}^M$ are mean free and therefore there are (unique) mean free $\hat{f}^k$ and $\hat{f}^M$ such that 
$$
\displaystyle
\LLs[\hat{f}^k]\ =\ \ks\,\hat{v}^k
\qquad\textrm{and}\qquad
\LLs[\hat{f}^M]=\ks\,\hat{v}^M\ .
$$
Then we denote $M_1$ the constant (with respect to $y$) so that 
\begin{equation}\label{sp_pb21}
\displaystyle
\tilde{\lambda}_1k_0+\tilde{\lambda}_0k_1+\ks\textrm{d}\cs(k_1,M_1)
=-6\ks\left(k_0<\hat{f}^k>+M_0<\hat{f}^M>\right).
\end{equation}
Substituting \eqref{bar_u1} into equation \eqref{eq_u2} and using the compatibility condition \eqref{cond_20} and relations \eqref{dku} yield then
\begin{equation}\label{def_u2}
\displaystyle
\hat{u}_2\ =\ \textrm{d}\Us(k_1,M_1)+k_0(<\hat{f}^k>-\hat{f}^k)+M_0(<\hat{f}^M>-\hat{f}^M)+\frac{k_2}{\ks}\Us'
\end{equation}
for some $k_2$ independent of $y$. Notice again $M_1=<\hat{u}_2>$. To complete the system on $(k_1,M_1)$, we still have to consider $\mathcal{O}(\nu^3)$ terms in \eqref{gKS_ls}
{\setlength\arraycolsep{1pt}
\begin{eqnarray}
\displaystyle
\LLs[\hat{u}_3]&+&\lambda_0\hat{u}_2+\ks\left((6\Us-\cs)\hat{u}_2+3\ks^2\frac{d^2\hat{u}_2}{dx^2}+\delta\left(2\ks\frac{d\hat{u}_2}{dx}+4\ks^3\frac{d^3\hat{u}_2}{dx^3}\right)\right)\nonumber\\
\label{nu3}
\displaystyle
&+&\lambda_1\hat{u}_1+\ks^2\left(3\ks\frac{d\hat{u}_1}{dx}+\delta\left(\hat{u}_1+6\ks^2\,\frac{d^2\hat{u}_1}{dx^2}\right)\right)\\
\displaystyle
&+&\lambda_2\hat{u}_0+\ks^3\hat{u}_0+4\delta\ks^4\frac{d\hat{u}_0}{dx}\ =\ 0\ .
\nonumber
\end{eqnarray}}
Setting $\tilde{\lambda}_1=\lambda_1/\ks$, the solvability condition for \eqref{nu3} expresses as
\begin{equation}\label{condu2}
\displaystyle
(\tilde{\lambda}_0-\cs)<\hat{u}_2>+\tilde{\lambda}_1<\hat{u}_1>+<6\Us\hat{u}_2>+\delta\ks<\hat{u}_1>\ =\ 0\ .
\end{equation}
Substituting \eqref{bar_u1} and \eqref{def_u2} into \eqref{condu2} yields
{\setlength\arraycolsep{1pt}
\begin{eqnarray}
\label{sp_pb22}
\displaystyle
(\tilde{\lambda}_0-\cs)\,M_1&+&\tilde{\lambda}_1M_0
+3\,\textrm{d}[<U^2>]_\star(k_1,M_1)\ =\ -\delta\ks M_0\\
\displaystyle
&+&<6\Us(\hat{f}^k-<\hat{f}^k>)>k_0+<6\Us(\hat{f}^M-<\hat{f}^M>)>M_0\ .
\nonumber
\end{eqnarray}}

It is easily seen by induction that we can obtain similary a complete set of equations for $(\tilde{\lambda}_n, k_n, M_n=<\hat{u}_{n+1}>)$ for any $n\in\N$ and thus obtain an expansion of both eigenvalues \emph{and} eigenvectors up to \emph{any order} in the small wavenumber limit $\nu\to 0$. We relate the two mode problems (\ref{cond_10},\ref{cond_20}) and (\ref{sp_pb21},\ref{sp_pb22}) with respectively "inviscid" and "viscous" Whitham's equations in the following section.

\subsubsection{Spectral validity of Whitham's equations}

In this section, we show the connection between the spectral analysis carried out previously and the spectral analysis of linearized Whitham's equations. Let us first check the validity of "inviscid" Whitham's equations,
\begin{equation}\label{wh_verif}
\displaystyle
\partial_T k+\partial_X\left(k\,c(k,M)\right)\ =\ 0,\quad
\partial_t M+\partial_X\left(<3U^2>\right)\ =\ 0\,.
\end{equation}
Linearizing (\ref{wh_verif}) about $(\ks,\Ms)$ in a co-moving frame $(X-\cs T,T)$ yields
\begin{equation}\label{wh_lin1}
\left\{
\begin{array}{rcl}
\displaystyle
\partial_T\tilde{k}&+&\ks\textrm{d}\cs\,(\partial_X\tilde{k},\partial_X\tilde{M})\ =\ 0\\[1ex]
\displaystyle
\partial_T\tilde{M}&+&\textrm{d}[<3U^2>]_\star(\partial_X\tilde{k},\partial_X\tilde{M})-\cs\partial_X\tilde{M}\ =\ 0
\end{array}
\right.\ .
\end{equation}
Performing a Fourier analysis, we look for solutions to equation \eqref{wh_lin1} as $(\tilde{k},\tilde{M})=e^{\nu(X+\Lambda T)}(\hat{\tilde{k}},\hat{\tilde{M}})$ and exactly obtain the linear system (\ref{cond_10},\ref{cond_20}) with $\Lambda=\tilde{\lambda}_0$ and $(\hat{\tilde{k}},\hat{\tilde{M}})=(k_0,M_0)$. As a result is also obtained the leading term in the expansion of the spectral curves associated to the stability of roll-waves, $\displaystyle \lambda(\nu)=\nu\ks\Lambda+\mathcal{O}(\nu^2)$. Thereby a necessary condition for roll-waves to be stable is that the Whitham's system \eqref{wh_verif} is \emph{hyperbolic}. In order to be able to check this hyperbolicity in the limit $\delta\to0$, we expand \eqref{disp_0}. Recall from \eqref{epMparam} that $c=c^{(0)}+6M$ and $<3U^2>=\bq+cM=\bq^{(0)}+3M^2$ with $c^{(0)}$ and $\bq^{(0)}$ independents of $M$. Thus the dispersion relation is
$$
\displaystyle
\Lambda^2+(\ks c^{(0)}{}'(\ks)-c^{(0)}(\ks))\Lambda-\left(\ks c^{(0)}(\ks)c^{(0)}{}'(\ks)+6\ks\bq^{(0)}{}'(\ks)\right)
\ =\ 0
$$
and system \eqref{wh_verif} is hyperbolic if
$$
\displaystyle
\Delta(\ks):=\left(c^{(0)}(\ks)+\ks c^{(0)}{}'(\ks)\right)^2+24\ks\,\bqs^{(0)}{}'(\ks)\ >\ 0\ .
$$ 

We further carry out the spectral analysis of the second order Whitham's equations. The first of these two equations reads
\begin{equation}\label{wh_v_verif}
\displaystyle
\partial_T k+\partial_X(k\,c(k,M))\ =\ -\varepsilon\partial_X[<k\,f>]
\end{equation}
where $f(k,M,\partial_Tk,\partial_Xk,\partial_TM,\partial_XM)$ is defined by $<f,U'>=0$ and 
{\setlength\arraycolsep{1pt}
\begin{eqnarray}
\displaystyle
\LL_{k,M}[f]&=&\partial_T (U-<U>)+3\,\partial_X\left(U^2-<U^2>\right)\nonumber\\
\displaystyle
&+&3\left(k^2\partial^3_{Xyy}U+k\partial_Xk\,\partial_y^2U\right)+\delta\left(2k\partial_{Xy}U+\partial_{X}k\,\partial_yU\right)\nonumber\\
&+&\delta\left(4k^3\partial^4_{Xyyy}U+6k^2\partial_X k\,\partial_y^3U\right)\ .
\nonumber
\end{eqnarray}}
Notice that the $\varepsilon$ dependence disappears, $\varepsilon$ being rescaled to $1$, by coming back to original variables $(x,t)=(X/\varepsilon,T/\varepsilon)$. Obviously one could have equivalently considered $\nu=\mathcal{O}(\varepsilon)$ wavenumbers in the Fourier plane wave analysis of linearized equations. Afterwards, linearizing \eqref{wh_v_verif} about $(\ks,\Ms)$ yields in a co-moving frame equation
\begin{equation}\label{lin_wh_v}
\displaystyle
\partial_t\tilde{k}+\ks\textrm{d}\cs(\partial_x\tilde{k},\partial_x\tilde{M})\ =\ -6\ks\,\partial_x<g>
\end{equation}
where $g=\textrm{d}f_{\star}(\tilde{k},\tilde{M},\partial_t\tilde{k}-\cs\partial_x\tilde{k},\partial_x\tilde{k},\partial_t\tilde{M}-\cs\partial_x\tilde{M},\partial_x\tilde{M})$ satisfies both conditions $<g,\,\Us'>=0$ and 
\begin{eqnarray}
\displaystyle
\LLs[g]&=&\partial_k\Us(\partial_t\tilde{k}-\cs\partial_x\tilde{k})
+\partial_k\left(3U^2-<3U^2>\right)_\star\partial_x\tilde{k}
\nonumber\\
\displaystyle
&&+\partial_M\left(3U^2-<3U^2>\right)_\star\partial_x\tilde{M}+3\left(\ks^2(\partial_k\Us)''+\ks\Us''\right)\partial_x\tilde{k}
\nonumber\\
\displaystyle
&&+\delta\left(\Us'+2\ks(\partial_k\Us)'+6\ks^2\Us'''+4\ks^3(\partial_k\Us)'''\right)\partial_x\tilde{k}
\nonumber\ .
\end{eqnarray}
We search for plane waves solution in the form $e^{\hat{\lambda}(\nu)t+\nu\ks x}\left(\hat{k}(\nu),\hat{M}(\nu)\right)$ and expand $\hat{\lambda}$, $\hat{k}$, $\hat{M}$ about $\nu=0$ as
$$
\begin{array}{rclclcl}
\displaystyle
\lambda(\nu)&=&\nu\ks \tilde{\lambda}_0&+&(\nu\ks)^2\tilde{\lambda}_1&+&\mathcal{O}(\nu^3)\ ,\\
\displaystyle
k(\nu)&=&k_0&+&\nu\ks k_1&+&\mathcal{O}(\nu^2)\ ,\\
M(\nu)&=&M_0&+&\nu\ks M_1&+&\mathcal{O}(\nu^2)\ .
\end{array}
$$
Identifying $\mathcal{O}(\nu)$ terms in \eqref{lin_wh_v} yields condition
$$
\displaystyle
\tilde{\lambda}_0 k_0+\ks\textrm{d}\cs(k_0,M_0)\ =\ 0
$$
which is again equation \eqref{cond_10}. We further identify $\mathcal{O}(\nu^2)$ terms and find
$$
\displaystyle
\tilde{\lambda}_1 k_0+\tilde{\lambda}_0 k_1+\ks\textrm{d}\cs(k_1,M_1)\ =\ -6\ks\left(<\hat{f}^M>k_0+<\hat{f}^M>M_0\right)
$$
which is \eqref{sp_pb21}. We focus now on the second Whitham's equation
\begin{equation}\label{wh_verif2}
\displaystyle
\partial_T M+\partial_X<3U^2>\ =\ \varepsilon\partial_X<6U(f-<f>)>-\varepsilon\delta\partial_{XX}M\ .
\end{equation}
Rescaling \eqref{wh_verif2} and linearizing about $(\ks,\Ms)$ in a co-moving frame yield
{\setlength\arraycolsep{1pt}
\begin{eqnarray}
\label{wh_lin2}
\displaystyle
\partial_t\tilde{M}&+&\textrm{d}[<3U^2>]_\star(\partial_X\tilde{k},\partial_X\tilde{M})-\cs\partial_X\tilde{M}\\
\displaystyle
&=&\partial_X<6\Us(\,g-<g>)>-\delta\partial_{XX}\tilde{M}\ .
\nonumber
\end{eqnarray}}
Identifying $\mathcal{O}(\nu)$ and $\mathcal{O}(\nu^2)$ terms yields respectively \eqref{cond_20} and \eqref{sp_pb22}.

$ $

As $\delta\to 0$, the domain of validity of the Whitham's equations shrinks to $0$ and one has to consider modulation with wavenumbers whose magnitudes are of order $\delta$. We analyze this situation in the next section.

\mathversion{bold}
\section{\label{sec4} Stability of wave trains for $\delta\to0$}
\mathversion{normal}

In this section, we consider the stability of periodic wave trains as $\delta\to 0$. In what follows, we will write $\delta=\varepsilon\bd$ where $\varepsilon$ is the characteristic modulation wavenumber. We first derive a set of three balance laws for modulated wave trains. In order to obtain a consistant spectral stability theory in the small wavenumber regime, and thus to be able to compare the Whitham's equations derived here and the spectral curves, we will have to consider wavenumbers $\nu=\mathcal{O}(\varepsilon)$. 

\subsection{Whitham's modulation equations}

Following section \ref{sec3}, we rescale (gKS) so that it reads
\begin{equation}\label{gKS4}
\displaystyle
\partial_T u+6u\partial_X u+\varepsilon^2\partial_X^3 u+\bd\left(\varepsilon^2\partial_X^2 u+\varepsilon^4\partial_X^4 u\right)\ =\ 0\ .
\end{equation}
We then search for modulated solutions in the form
\begin{eqnarray*}
\displaystyle
u(X,T)&=&u_0\left(\frac{\varphi(X,T)}{\varepsilon};X,T\right)+\varepsilon u_1\left(\frac{\varphi(X,T)}{\varepsilon};X,T\right)+\mathcal{O}(\varepsilon^2)\ ,\\
\displaystyle
\varphi(X,T)&=&\varphi_0(X,T)+\varepsilon \varphi_1(X,T)+\mathcal{O}(\varepsilon^2)\ ,
\end{eqnarray*}
with $u_0(y;X,T)$ and $u_1(y;X,T)$ $1$-periodic in $y$. Inserting this ansatz into \eqref{gKS4} and identifying $\mathcal{O}(\varepsilon^{-1})$ yield after setting $k_0=\partial_x\varphi_0$ and $\Omega_0=\partial_T\varphi_0$
\begin{equation}\label{sl-1}
\displaystyle
\Omega_0\partial_y u_0+6k_0 u_0\partial_y u_0+k_0^3\partial_y^3 u_0\ =\ 0\ .
\end{equation}
Thus $u_0(\,\cdot\,;X,T)$ is a periodic traveling wave of (KdV). These solutions are given in term of elliptic functions (see section \ref{sec2}),
\begin{eqnarray*}
\displaystyle
u_0(y;X,T)&=&M_0(X,T)+[\tilde{U}_0(k_0(X,T),p_0(X,T))](y)\ ,\\
c_0(X,T)&=&6M_0(X,T)+\tilde{c}_0(k_0(X,T),p_0(X,T))\ ,
\end{eqnarray*}
for some $(M_0,p_0)$, with $\tilde{U}_0$ and $\tilde{c}_0$ given by \eqref{per_kdv}. From now on, we drop the $0$ subscript in the definitions \eqref{per_kdv} of $\tilde{U}$ and $\tilde{c}$ and keep it for the functions of paremeters $(k,M,p)$ evaluated in $(k_0,M_0,p_0)$. Moreover we will denote\footnote{This notation is of course not consistant with the meaning of $U$ and $c$ up to now.} $U(k,M,p)=M+\tilde{U}(k,p)$ and $c(k,M,p)=6M+\tilde{c}(k,p)$. Compatibility condition 
$\displaystyle \partial_T\partial_X\varphi_0=\partial_X\partial_T\varphi_0$ yields then
\begin{equation}\label{w4_0}
\displaystyle
\partial_T k_0+\partial_X\left(k_0\,c(k_0,M_0,p_0)\right)\ =\ 0\ .
\end{equation}
Next, we consider $\mathcal{O}(\varepsilon^0)$ terms to obtain a system for $(k_0,M_0,p_0)$ in a closed form. Setting $\Omega_1=\partial_T\varphi_1$ and $k_1=\partial_X\phi_1$ yields
{\setlength\arraycolsep{1pt}
\begin{eqnarray}
\displaystyle
\LLkdv_0[u_1]&+&\Omega_1\partial_y u_0+k_1\left(3k_0^2\partial_y^3u_0+6u_0\partial_y u_0\right)+\partial_T u_0+\partial_X(3u_0^2)\nonumber\\
\label{sl0}
\displaystyle
&+&3(k_0^2\partial_{Xyy}u_0+k_0\partial_X k_0\,\partial_{y}^2u_0)+\bd\left(k_0^2\partial_y^2u_0+k_0^4\partial_y^4u_0\right)=0
\end{eqnarray}}
with $\displaystyle \LLkdv_{k,p}[f]=k\partial_y\left((6\tilde{U}(k,p)-\tilde{c}(k,p))f\right)+k^3\partial_y^3f$. Equation \eqref{sl0} has a solution if and only two compatibility conditions are satisfied. Indeed $\LLkdv$ is Fredholm of index $0$ and the kernel of $\LLkdv^{\textrm{ad}}$ is spanned by $1$ and $U$ (or $\tilde{U}$). The first compatibility condition yields
\begin{equation*}
\displaystyle
\partial_T <u_0>+\partial_X <3u_0^2>\ =\ 0
\end{equation*}
thus
\begin{equation}\label{w4_1}
\displaystyle
\partial_T M_0\,+\,3\,\partial_X<U_0^2>\ =\ 0\ .
\end{equation}
The second compatibility condition reads
{\setlength\arraycolsep{1pt}
\begin{eqnarray*}
\displaystyle
\!\!\!\!\!
<\partial_T u_0,U_0>&+&<\partial_X(3u_0^2),U_0>+3<k_0^2\partial_{Xyy}u_0+k_0\partial_X k_0\,\partial_y^2 u_0,U_0>\\
\displaystyle
&+&\bd\left(k_0^2<\partial_y^2 u_0,U_0>+k_0^4<\partial_y^4 u_0,U_0>\right)=0.
\end{eqnarray*}
where $<\,\cdot\,,\,\cdot\,>$ is the usual scalar product on $1$-periodic functions. Integrating by parts yields equation
\begin{equation}
\displaystyle\label{w4_2}
\frac12\partial_T<U_0^2>\,+\,\partial_X <2U_0^3-\frac32k_0^2(U_0')^2>
\,=\,\bd k_0^2<(U_0')^2-k_0^2(U_0'')^2>.
\end{equation}

As a result, one obtains one \emph{extra balance law} in comparison to the case where $\delta>0$ is fixed. Indeed, in the case $\delta\to 0$, (gKS) is a perturbation of (KdV) and for $\bd=0$, one exactly obtains the modulated system for (KdV), as it was obtained in \cite{Z_kdv}. However, in their derivation, they considered the energy equation
$$
\displaystyle 
\partial_t\left(\frac{u^2}{2}\right)+\partial_x\left(2u^3+u\partial_x^2 u-\frac{(\partial_x u)^2}{2}\right)\ =\ 0
$$
obtained from (KdV) by multiplying the equation by $u$, as an extra conservation law and they averaged this conservation law to obtain \eqref{w4_2} (when $\bd=0$). Notice in particular that many conservation laws were available for such a process. Here, this extra equation seems to be obtained more naturally as a comptability condition in the process of computing asymptotic expansions of low frequency perturbations of the wave. The stationnary solutions to (\ref{w4_0}, \ref{w4_1},\ref{w4_2}) are exactly the parameters $k_0,M_0,p_0$ such that the following relation holds
\begin{equation}\label{w4_r}
\displaystyle
<\left(\tilde{U}(k_0,p_0)'\right)^2>\quad=\quad k_0^2<\left(\tilde{U}(k_0,p_0)''\right)^2>\ .
\end{equation}
It corresponds exactly to the solvability condition in order to compute an expansion of periodic solution to (gKS) as $\delta\to 0$. This restricts the set of periodic solutions to a two dimensional manifold just as in the case $\delta>0$ fixed. From the dynamical point of view, we recover a set of two equations by considering the relaxation limit $\bd\to\infty$. This yields the relaxed system (\ref{w4_0},\ref{w4_1},\ref{w4_r}) which is a system of two conservation laws, just as in the case of fixed  $\delta>0$.

Should conditions \eqref{w4_1} and \eqref{w4_2} hold, one may solve equation
\begin{eqnarray*}
\displaystyle
\LLkdv_0[v_0]&+&\partial_T u_0+\partial_X(3u_0^2)\\
\displaystyle
&+&3(k_0^2\partial_{Xyy}u_0+k_0\partial_X k_0\,\partial_{y}^2u_0)+\bd\left(k_0^2\partial_y^2u_0+k_0^4\partial_y^4u_0\right)=0
\end{eqnarray*}
with $v_0$, in a unique way when added conditions\footnote{See the description of the kernel of $\LLkdv$ in the next subsection.}
$$\displaystyle<v_0,U'_0>\ =\ <v_0>\ =\ 0\ .$$ 
Note that the fonction $v_0$ is odd\footnote{Of course most of unobvious claims of oddness or evenness may be deduced from uniqueness.}. Equation \eqref{sl0} reads then
$$
\displaystyle
\LLkdv_0[u_1-v_0-k_1\partial_k U_0]+(\Omega_1-k_1\partial_k \Omega_0)\,U'_0\ =\ 0\ .
$$
We denote $M_1=<u_1>$ and introduce $p_1$ the constant (with respect to $y$) such that  $\displaystyle\Omega_1=\textrm{d}\Omega_0[k_1,M_1,p_1]$. Solutions to \eqref{sl0} are then written as
$$
\displaystyle
u_1=\textrm{d}U_0[k_1,M_1,p_1]+v_0+\gamma_1 U'_0
$$
for some $\gamma_1$ independent of $y$. In what follows, we denote $u_1^e=\textrm{d}U_0[k_1,M_1,p_1]$ and $u_1^o=v_0+\gamma_1 U'_0$, respectively even and odd parts of $u_1$. Compatability condition $\partial_T\partial_X\phi_1=\partial_X\partial_T\phi_1$ yields 
\begin{equation}\label{evol_k1}
\displaystyle
\partial_T k_1-\partial_X\left(d\Omega_0[k_1,M_1,p_1]\right)\ =\ 0\ .
\end{equation}
In order to obtain the full set of equations that governs the evolution of $(k_1,p_1,M_1)$, one has to consider $\mathcal{O}(\varepsilon)$ terms in \eqref{gKS4} and finds
\begin{eqnarray}
\displaystyle
\LLkdv_0[u_2]&+&\Omega_1\partial_yu_1^o+\partial_T u_1^e+\partial_X(6u_0u_1^e)\nonumber\\
\displaystyle
&+&6\partial_y(k_1u_0 u_1^e+k_0u_1^eu_1^o)+3(k_0^2\partial_{Xyy}u_1^e+k_0\partial_Xk_0\,\partial_y^2u_1^e)\nonumber\\
\label{sl1}
\displaystyle
&+&3\left(2k_0k_1\partial_{Xyy}u_0+\partial_X(k_0k_1)\partial_{yy}u_0+k_0^2k_1\partial_y^3 u_1^o\right)\\
\displaystyle
&+&\bd\left(k_0^2\partial_y^2u_1^e+k_0^4\partial_y^4u_1^e+2k_0k_1\partial_y^2u_0+4k_0^3k_1\partial_y^4u_0\right)+F^o=0
\nonumber
\end{eqnarray}
where $F^o$ is some odd function. First solvability condition for \eqref{sl1} reads
$$
\displaystyle
\partial_T<u_1^e>\ +\ \partial_X<6u_0u_1^e>\ =\ 0\ .
$$
Inserting $u_1^e$ into this latter equation yields
\begin{equation}\label{evol_M1}
\displaystyle
\partial_T M_1+\partial_X\left(\textrm{d}[<3U^2>]_0[k_1,M_1,p_1]\right)\ =\ 0\ .
\end{equation}
The second solvability condition gives
{\setlength\arraycolsep{1pt}
\begin{eqnarray}
\displaystyle
<\partial_Tu_1^e,U_0>&+&<\partial_X(6u_0u_1^e),U_0>-<6k_1u_0u_1^o+6k_0u_1^eu_1^o+\Omega_1 u_1^o,U_0'>\nonumber\\
\displaystyle
&-&\frac32\left(2<k_0^2\partial_{Xy}u_1^e,U_0'>+\,\partial_X\left(k_0^2\right)<\partial_{y}u_1^e,U_0'>\right)\nonumber\\
\label{evol_p1}
\displaystyle
&-&3\left(\partial_X<k_0k_1(U_0')^2>-k_0^2k_1<\partial_y u_1^o,U_0''>\right)\\
\displaystyle
&-&\bd\left(k_0^2<\partial_y u_1^e,U_0'>-k_0^4<\partial_y^2 u_1^e,U_0''>\right)\nonumber\\
&-&\bd\left(k_1\left(2k_0<(U_0')^2>-4k_0^3<(U_0'')^2>\right)\right)\ =\ 0\ .
\nonumber
\end{eqnarray}}
Equations (\ref{evol_k1},\ref{evol_M1},\ref{evol_p1}) govern the evolution of $(k_1,M_1,p_1)$.

If these equations are satisfied, we can further solve \eqref{sl1} and express $u_2$ in terms of $(k_2,M_2,p_2)$ and $(k_i,p_i,M_i,u_i)_{i=0,1}$. Actually, by induction, one can proceed with expansions up to \emph{any order} with respect to $\varepsilon$.

Unlike the case $\delta>0$ fixed, by summing (\ref{w4_0},\ref{w4_1},\ref{w4_2}) and $\varepsilon$(\ref{evol_k1},\ref{evol_M1},\ref{evol_p1}) we can not obtain a closed system governing the evolution of $(k,M,p)=(k_0,M_0,p_0)+\varepsilon(k_1,M_1,p_1)$ up to order $\mathcal{O}(\varepsilon^2)$ terms. The first two equations would indeed close into
\begin{equation*}
\begin{array}{rcl}
\displaystyle
\partial_T k-\partial_X\Omega&=&0\ ,\\
\displaystyle
\partial_T M+\partial_X<3U^2>&=&0\ ,
\end{array}
\end{equation*}
but the last equation would not. However, we will relate the stability of periodic wave trains in the limit $\delta\to 0$ with a wave analysis of the two first steps of the modulation.
 
\subsection{Bloch analysis of spectral stability}
 
In this section, we analyse the spectral stability of periodic wave trains in the limit $\delta\to 0$. Recall that they expand as 
$$
{\setlength\arraycolsep{0.5ex}
\begin{array}{rclclclclcl}
\displaystyle
U^{\delta}(y;k,M)&=&M&+&\tilde{U}(y;k,p)&+&\delta\tilde{U}_1(y;k,p)&+&\delta^2\tilde{U}_2(y;k,p)&+&\mathcal{O}(\delta^3)\\[1ex]
\displaystyle
c^{\delta}(k,M)&=&6M&+&\tilde{c}\,(k,p)&+&\delta^2\tilde{c}_2(k,p)&+&\mathcal{O}(\delta^3)&&
\end{array}}
$$
with parameter $p$ obtained from $k$ through equation \eqref{eq9}. Moreover stand the following relations
\begin{equation}\label{dl12}
\begin{array}{l}
\displaystyle
k^3\tilde{U}_1'''+k\left((6\tilde{U}-\tilde{c})\tilde{U}_1\right)'+k^2\tilde{U}''+k^4\tilde{U}''''\ =\ 0\ ,\\
\displaystyle
k^3\tilde{U}_2'''+k\left((6\tilde{U}-\tilde{c})\tilde{U}_2+3\tilde{U}_1^2-\tilde{c}_2\tilde{U}\right)'+k^2\tilde{U}''_1+k^4\tilde{U}''''_1\ =\ 0\ .
\end{array}
\end{equation}
Notice $\tilde{U}_1$ is odd and $\tilde{U}_2$ is even.

The spectral problem associated to \eqref{gKS4} is written in a co-moving frame
\begin{equation}\label{gKS4_l}
\begin{array}{rcl}
\displaystyle
\lambda(\nu)\hat{u}&+&\ks\left(\frac{d}{dx}+\nu\right)\left(\left(6\Us^{\varepsilon\bd}-\cs^{\varepsilon\bd}\right)\hat{u}\right)
+\ks^3\left(\frac{d}{dx}+\nu\right)^3\hat{u}\\[1ex]
&+&\varepsilon\bd\left(\ks^2\left(\frac{d}{dx}+\nu\right)^2\hat{u}+\ks^4\left(\frac{d}{dx}+\nu\right)^4\hat{u}\right)\ =\ 0
\end{array}
\end{equation}
with a $1$-periodic $\hat{u}$. In the limit $\varepsilon\to0$ the analysis of \eqref{gKS4_l} has been carried out in \cite{Bar} both for $\nu=\mathcal{O}(1)$ and for $\nu=\mathcal{O}(\varepsilon)$ (in an appendix). Here we recall only the latter one since this is the one to be interpreted in terms of linearized modulation equations. Thus we denote $\nu=\bnu\varepsilon$ and expand eigenfunctions and eigenvalues as
$$
\displaystyle
\lambda(\nu)\ =\ \varepsilon\bnu\lambda_0+(\varepsilon\bnu)^2\lambda_1+\mathcal{O}(\varepsilon^3)\ ,\qquad
\hat{u}\ =\ \hat{u}_0+\varepsilon\hat{u}_1+\varepsilon^2\hat{u}_2+\mathcal{O}(\varepsilon^3)\ .
$$
Identifying $\mathcal{O}(\varepsilon^0)$ terms in \eqref{gKS4_l}, one finds
$$
\displaystyle
\LLkdvs[\hat{u}_0]\ =\ \ks\left((6\tilde{U}_\star-\tilde{c}_\star)\hat{u}_0\right)'+\ks^3\hat{u}_0'''\ =\ 0\ .
$$
The two dimensional kernel of $\LLkdvs$ is spanned by $U'_\star$ and 
$$
V_\star\ =\ \partial_M U_\star-\frac{\partial_M\cs}{\partial_p\cs}\partial_p U_\star\ =\ 1-6[\partial_p\cs]^{-1}\partial_p U_\star\ .
$$
Thus $\hat{u}_0=k_1 U_\star'+M_0 V_\star$ for some $(k_1,M_0)$. Notice $<\hat{u}_0>=M_0$. Next, we consider $\mathcal{O}(\varepsilon)$ terms
\begin{equation}\label{gks4l1}
\begin{array}{rcl}
\displaystyle
-\LLkdvs[\hat{u}_1]&=&\ks\bd\left(6\tilde{U}_{1\star}\,\hat{u}_0\right)'+\,\bd\left(\ks^2\hat{u}_0''+\ks^4\hat{u}_0''''\right)\\[1ex]
&&+\ 3\bnu\ks^3\hat{u}_0''\ +\ \bnu\ks\left(\tilde{\lambda}_0\hat{u}_0+(6\tilde{U}_\star-\tilde{c}_\star)\hat{u}_0\right)
\end{array}
\end{equation}
where again $\tilde{\lambda}_0$ stands for $\lambda_0/\ks$. Solvability conditions are $<$\eqref{gks4l1}${}_{\textrm{r.h.s.}}>=0$ and $<U_\star,$\eqref{gks4l1}${}_{\textrm{r.h.s.}}>=0$. They reduce to $M_0=0$. Equation \eqref{gks4l1} then reads
$$
\displaystyle
\LLkdvs[\hat{u}_1-\bnu\ks k_1\partial_k U_\star-k_1\bd\tilde{U}_{1,\star}']+\bnu\ks k_1\left(\tilde{\lambda}_0+\ks\partial_k\cs\right)U_\star'\ =\ 0\ .
$$
We introduce $M_1=<\hat{u}_1>$ then $p_1$ the constant such that
\begin{equation}\label{lin1}
\displaystyle
\tilde{\lambda}_0\,k_1+\ks\textrm{d}\cs(k_1,M_1,p_1)\ =\ 0
\end{equation}
where $c(k,M,p)=6M+\tilde{c}(k,p)$. Therefore accordingly
\begin{equation}\label{hat_u1}
\displaystyle
\hat{u}_1\ =\ \bnu\ks\,\textrm{d}U_\star(k_1,M_1,p_1)+k_1\bd\tilde{U}_{1,\star}'+k_2 U_\star'
\end{equation}
where again $U(k,M,p)=M+\tilde{U}(k,p)$. To determine $(\tilde{\lambda}_0,k_1,M_1,p_1)$, one must look at \eqref{gKS4_l} up to order $\mathcal{O}(\varepsilon)$ and then finds
\begin{equation}
\label{gks4l2}
\begin{array}{rclcl}
\displaystyle
\!\!\!\!
-\LLkdvs[\hat{u}_2]&=&\bnu(\lambda_0\hat{u}_1+\bnu\lambda_1\hat{u}_0)
&+&\ks\bd\left(6\tilde{U}_{1,\star}\hat{u}_1+\bd(6\tilde{U}_{2,\star}-\tilde{c}_{2,\star})\hat{u}_0\right)'
\\
\displaystyle
&+&3\ks^3\bnu\left(\hat{u}_1''+\bnu\hat{u}_0'\right)
&+&\bnu\ks\left((6U_\star-\cs)\hat{u}_1+6\bd\,\tilde{U}_{1,\star}\hat{u}_0\right)
\\
\displaystyle
&+&\bd\ks^2(\hat{u}_1''+2\bnu\hat{u}_0')
&+&\bd\ks^4(\hat{u}_1''''+4\bnu\hat{u}_0''')\ .
\end{array}
\end{equation}
Solvability conditions are $<$\eqref{gks4l2}${}_{\textrm{r.h.s.}}>=0$ and $<U_\star,$\eqref{gks4l2}${}_{\textrm{r.h.s.}}>=0$.
Taking into account \eqref{hat_u1} first condition reads
\begin{equation}\label{lin2}
\displaystyle
(\tilde{\lambda}_0-\cs)\, M_1\ +\ \textrm{d}[<3U^2>]_\star(k_1,M_1,p_1)\ =\ 0\ .
\end{equation}
With odd and even parts notations, equation \eqref{gks4l2} is also written
\begin{equation*}
\begin{array}{rcl}
\displaystyle
\!\!\!\!
-\LLkdvs[\hat{u}_2]&=&\bnu\lambda_0\hat{u}_1^e
+\ks\bd\left(6\tilde{U}_{1,\star}\hat{u}_1^e+\bd(6\tilde{U}_{2,\star}-\tilde{c}_{2,\star})\hat{u}_0^o\right)'
\\
\displaystyle
&+&3\ks^3\bnu\left(\hat{u}_1^e{}''+\bnu\hat{u}_0^o{}'\right)
+\bnu\ks\left((6U_\star-\cs)\hat{u}_1^e+6\bd\,\tilde{U}_{1,\star}\hat{u}_0^o\right)
\\
\displaystyle
&+&\bd\ks^2(\hat{u}_1^e{}''+2\bnu\hat{u}_0^o{}')
+\bd\ks^4(\hat{u}_1^e{}''''+4\bnu\hat{u}_0^o{}''')\ +\ L^o\ .
\end{array}
\end{equation*}
where $L^o$ is some odd term. Using the second equation of \eqref{dl12} and the fact that
$$
<U_\star,\ks\tilde{U}_{2,\star}''''+\ks\left((6U_\star-\cs)\tilde{U}_{2,\star}'\right)'>\ =\ 0
$$
(coming from $\LLkdvs U'_\star=0$), the last solvability condition is thus written
\begin{equation}
\label{lin3}
\begin{array}{rcl}
\displaystyle
\!\!\!\!
0&=&
\bnu\ks\left[(\tilde{\lambda}_0-\cs)\textrm{d}<\frac12U^2>_\star
+\textrm{d}<2U^3-\frac{3}{2}k^2(U')^2>_\star
\right](k_1,M_1,p_1)
\\[1ex]
\displaystyle
&-&\bd\ \textrm{d}<k^2(U')^2-k^4(U'')^2>_\star (k_1,M_1,p_1)
\\[1ex]
\displaystyle
&-&6\ks\bd\ <\tilde{U}_{1,\star}U'_\star\textrm{d}U_\star(k_1,M_1,p_1)>
\\[1ex]
\displaystyle
&+&k_1\bd\left((\tilde{\lambda}_0-\cs)<\tilde{U}_{1,\star}'U_\star>-<3\ks^2\tilde{U}_{1,\star}''+6\tilde{U}_{1,\star}U_\star,U'_\star>\right)\ .
\end{array}
\end{equation}

As a result, (\ref{lin1},\ref{lin2},\ref{lin3}) forms a full spectral problem  for eigenvalue $\tilde{\lambda}_0$ and eigenvector $(k_1,M_1,p_1)$. Therefore $\tilde{\lambda}_0$ must satisfy a dispersion relation $E(\tilde{\lambda}_0)=0$ with $E$ a polynom of degree $3$. This dispersion is analogous to the one found in \cite{Bar}. In the following, we interpret the system (\ref{lin1},\ref{lin2},\ref{lin3}) with the help of a wave analysis of modulation systems derived in the previous section.
 
\subsection{Spectral validity of Whitham's equations}

Let us first remark that a byproduct of the previous analysis is the direct spectral validation of modulated Whitham's equations for (KdV). We recover (KdV) setting $\bd=0$ in previous computations. The linear spectral problem (\ref{lin1},\ref{lin2},\ref{lin3}) reads then
\begin{equation*}
\left\{
{\setlength\arraycolsep{1pt}
\begin{array}{rcl}
\displaystyle
\Big[(\tilde{\lambda}_0-\cs)\,\textrm{d}\ks
&+&\textrm{d}[kc]_\star\Big](k_0,M_0,p_0)\ =\ 0
\\
\Big[(\tilde{\lambda}_0-\cs)\,\textrm{d}\Ms
&+&\textrm{d}[<3U^2>]_\star\Big](k_0,M_0,p_0)\ =\ 0
\\
\displaystyle
\Big[\frac12(\tilde{\lambda}_0-\cs)\textrm{d}[<U^2>]_\star
\displaystyle
&+&\textrm{d}\big[<2U^3>-\frac32 k<(U')^2>\big]_\star\Big](k_0,M_0,p_0)\,=\,0
\end{array}}
\right.
\end{equation*}
whereas the modulated system for (KdV), obtained by setting $\bd=0$ in (\ref{w4_0},\ref{w4_1},\ref{w4_2}), is
\begin{equation}\label{kdv_m}
\left\{\begin{array}{rcl}
\displaystyle
\partial_T k_0&+&\partial_X(k_0\,c_0)=0\\[1ex]
\displaystyle
\partial_T M_0&+&\partial_X <3U_0^2>=0\\[1ex]
\displaystyle
\frac12\partial_T<U_0^2>
&+&\partial_X <2U_0^3-\frac32 k_0^2(U'_0)^2>=0
\end{array}\right.\ .
\end{equation}
There is no difficulty here to see that the spectral problem described above and the one obtain when considering the stability of steady solutions to \eqref{kdv_m} are the same. This validates modulated equations for (KdV) on a spectral level. Moreover, hyperbolicity of modulated equations \eqref{kdv_m} is thus a necessary condition for periodic travelling waves solutions to (KdV) to be stable. 

The situation is trickier for (gKS): from the spectral analysis carried out in the previous section, we see that the correction $\tilde{U}_1$ to $U$ plays a significant role. As a result, the linearization of (\ref{w4_0},\ref{w4_1},\ref{w4_2}) at a steady solution is clearly not sufficient and we need to take into account higher order modulations. Actually the right way to search for "solutions" in the form $(\ks,\Ms,\ps)+\varepsilon (k_1,M_1,p_1)$, with $(\ks,\Ms,\ps)$ constant and $(k_1,M_1,p_1)$ a plane wave, is to set in the modulation process $(k_0,M_0,p_0)=(\ks,\Ms,\ps)$ and look for $(k_1,M_1,p_1)$ as a plane wave. Obviously, when $\bd=0$, it coindides with the direct approach.

Setting $(k_0,M_0,p_0)=(\ks,\Ms,\ps)$ (with $\ps$ obtained through \eqref{eq9}), equations \eqref{evol_k1} and \eqref{evol_M1} turn into
\begin{eqnarray}
\label{mod_lin1}
\displaystyle
\partial_T k_1\ +\ \textrm{d}[kc]_\star\partial_X(k_1,M_1,p_1)&=&0\\
\label{mod_lin2}
\displaystyle
\partial_T M_1\ +\ \textrm{d}<3U^2>_\star\partial_X(k_1,M_1,p_1)&=&0
\end{eqnarray}
and, since then $v_0=\bd\tilde{U}_{1,\star}$, using 
$$
<U_\star,\LLkdvs (\textrm{d}U'_\star(k_1,M_1,p_1))>\ =\ 0\ ,
$$
condition \eqref{evol_p1} reads
\begin{equation}
\label{mod_lin3}
\begin{array}{rcl}
\displaystyle
\!\!\!\!
0&=&\textrm{d}<\frac12U^2>_\star\partial_T(k_1,M_1,p_1)
+\textrm{d}<2U^3-\frac{3}{2}k^2(U')^2>_\star\partial_X(k_1,M_1,p_1)
\\[1ex]
\displaystyle
&-&\bd\ \textrm{d}<k^2(U')^2-k^4(U'')^2>_\star\partial_X(k_1,M_1,p_1)
\\[1ex]
\displaystyle
&-&6\ks\bd\ <\tilde{U}_{1,\star}U'_\star\,\textrm{d}[kU]_\star(k_1,M_1,p_1)>
\\[1ex]
\displaystyle
&-&k_1\bd<3\ks^2\tilde{U}_{1,\star}'',U'_\star>-\bd<\tilde{U}_{1,\star}'U_\star>\textrm{d}[kc]_\star(k_1,M_1,p_1)\ .
\end{array}
\end{equation}

Searching, in a co-moving frame, for plane-wave solutions to (\ref{mod_lin1},\ref{mod_lin2},\ref{mod_lin3}) leads exactly to the spectral problem (\ref{lin1},\ref{lin2},\ref{lin3}). This validates at the spectral level the modulation equations in the limit $\delta\to 0$.

\section{Conclusion}

In this paper, we derived Whitham's equations for the modulations of periodic traveling waves, solutions to generalized Kuramoto Sivashinski (gKS) equations. We proved that this set of equations gives necessary conditions of first and second order for the spectral stability of travelling waves in the low frequency regime. Indeed, the hyperbolicity of Whitham's first order equations implies the tangency of spectral curves at the origin whereas the parabolicity of the second order Whitham's system ensures that the spectral curves lie in the stable spectral region. We also reinterpreted the spectral analysis carried out in \cite{Bar} in the (KdV) limit $\delta\to 0$ with the help of a set of Whitham's equations. As a byproduct, we have obtained a set Whitam's equations for modulated wavetrains in (KdV) equation in a more natural way than in  \cite{Z_kdv}  where an additional equation for the energy was introduced to close the Whitham's set of equations.\\

\noindent
It is worth noting that, following the approach already introduced in \cite{D3S} in the case of reaction diffusion equations, one could have also derived a viscous Burgers equation by choosing the scaling:
$$
\begin{array}{ll}
\displaystyle
u(X,T)=\sum_{k=0}^{\infty}\varepsilon^ku_k\big(\frac{\phi(X,T)}{\varepsilon},X-c_gT,\varepsilon T),\big),\\
\displaystyle
\phi(X,T)=\phi_\star(X,T)+\sum_{k\geq 1}\varepsilon^k \phi^k(X-c_gT,\varepsilon T),
\end{array}
$$
\noindent
where $\phi_\star$ is the phase associated to a fixed periodic wave train $U_{\star}(.,k_{\star},M_{\star})$ and $c_g$ is one of the characteristic speeds associated to the inviscid Whitham's system. This is the first step towards the construction of approximate solutions to (gKS) on asymptotically large time interval $O(T^{\star}/\varepsilon^2)$ with arbitrary $T^{\star}>0$. Furthermore, as viscous Burgers equations have shock structures, we expect that we will be able to construct solutions of shock type in (gKS) in the form of a travelling wave relating two families of periodic wave train at $-\infty$ and $+\infty$. This approach will be developped in a forthcoming paper.\\

\noindent
Another problem of interest is the stability of periodic travelling waves, solutions to (gKS): we expect that, using the approach introduced in \cite{Jo_Zu}, one can prove the nonlinear stability of periodic wavetrains solutions to (gKS) under suitable spectral assumptions. Concerning this spectral stability issue, there should be a range of periods for which periodic wavetrains are stable, just as in the case of roll-waves solutions to shallow water equations. Indeed, at the onset of periodic wave train (close to the Hopf bifurcation), steady solutions are unstable and by perturbation arguments, we shall prove the same result for small amplitude solutions. On the other hand, solitary waves are also unstable since, for (gKS), constants are unstable under low frequency perturbations. This heuristic approach was partially validated in \cite{Bar} in the limit $\delta\to 0$ but the authors only considered perturbations with zero mean. This analysis should be completed for general perturbations in order to conclude to the full stability of periodic wavetrains.


\end{document}